\documentclass[11pt]{article}
\usepackage[utf8]{inputenc}
\usepackage[tbtags]{amsmath}
\usepackage{amssymb}
\usepackage{amsthm}
\usepackage[misc]{ifsym}
\usepackage{cases}
\usepackage{mathrsfs}
\usepackage{color}
\usepackage{graphicx}
\numberwithin{equation}{section}   
\normalsize
\title{\bf Linear Quadratic  Control of Backward Stochastic Differential Equation with Partial Information \thanks{This work  supported by   the National Natural Science Foundations of China under Grants 61821004, 61633015, 61877062,  and 61977043.}}

\author{\normalsize Guangchen Wang\thanks{\it School of Control Science and Engineering, Shandong University, Jinan 250061, PR China, E-mail: wguangchen@sdu.edu.cn}, Wencan Wang\thanks{\it School of Control Science and Engineering, Shandong University, Jinan 250061, PR China, E-mail: wwencan@163.com}, Zhiguo Yan\thanks{\it School of Electrical Engineering and Automation, Qilu University of Technology (Shandong Academy of Sciences), Jinan 250353, PR China, E-mail:yanzg500@sina.com}}
\newtheorem{theorem}{Theorem}[section]
\newtheorem{lemma}{Lemma}[section]

\newtheorem{remark}{Remark}[section]
\begin{document}
\maketitle

\noindent{\bf Abstract:}\quad In this paper, we study  an    optimal control problem  of linear backward stochastic differential equation (BSDE) with quadratic  cost functional under partial information.  This problem is solved completely and explicitly by using a stochastic maximum principle and a decoupling technique. By using the  maximum principle, a stochastic Hamiltonian system, which is a forward-backward stochastic differential equation (FBSDE) with filtering, is obtained. By decoupling the stochastic Hamiltonian system, three Riccati equations, a  BSDE with filtering, and a    stochastic differential equation (SDE) with filtering are derived. We then get an optimal control with a  feedback representation. An explicit  formula  for the  corresponding optimal cost is also established. As illustrative examples, we consider two   special scalar-valued control problems and give some numerical simulations.

\vspace{2mm}

\noindent{\bf Keywords:}\quad Linear quadratic optimal control; backward stochastic differential equation; filtering; Ricatti equation; feedback representation.

\vspace{2mm}

\noindent{\bf Mathematics Subject Classification:}\quad 93E20, 60H10

\section{Introduction}

A  BSDE is an It\^o SDE for which a random terminal rather than an initial condition on
state has been specified. Bismut \cite{J2000An} first introduced a linear BSDE, which is an  adjoint equation  of  stochastic
optimal control problem. Pardoux and Peng \cite{Pardoux1990Adapted} extended the linear BSDE to a general case. Since then, there has been considerable attention    on related topics and their applications  among researchers in mathematical finance and stochastic optimal control. See for example, El Karoui et al. \cite{Karoui1997Backward}, Ma and Yong \cite{Ma2007Forward}, Kohlmann and Zhou \cite{Kohlmann2000Relationship}.

Since BSDE  stems from stochastic control theory, it is very natural and appealing to investigate  optimal  control  problem of BSDE. Moreover,   controlled BSDE  is expected to   have wide and important applications in various fields, especially in mathematical finance.  In financial investment,  a European contingent claim  $\xi$, which is a random variable, can be thought as a contract to be guaranteed at maturity $T$. Peng \cite{Peng1993Backward} and  Dokuchaev and Zhou \cite{Nikolai1999Stochastic} derived  local and global stochastic maximum principles of optimality for BSDEs, respectively. Linear quadratic (LQ)  optimal control problems of BSDEs have also been investigated.  Lim and Zhou \cite{Lim2001Linear} discussed an LQ control problem of BSDE with a general setting  and gave a feedback representation of the optimal control. Li et al. \cite{Li2017Linear} extended the results in \cite{Lim2001Linear} to the case with mean-field term. Huang et al. \cite{Huang2016Backward} and Du et al. \cite{Kai2018Linear} considered LQ backward mean-field games. Du and Wu \cite{Du2019Linear} concerned  a    stackelberg game for mean field linear BSDE with   quadratic cost functionals.   

In this paper, we investigate an  LQ   control problem  of BSDE with partial information, which will be referred as a stochastic backward LQ control problem. We are devoted to deriving the optimal control with a   feedback representation and establishing an explicit formula for the corresponding optimal cost. Note that the mentioned papers above are concentrated on the complete information case. The motivation of studying
stochastic  control problems  with partial information arises naturally from the area of financial  economics. In a  portfolio and  consumption problem, let $\mathbb F\equiv\{\mathcal F_t\}_{t\geq0}$  denote the flow of information  generated by all market noises. In reality,  the  information available to  an agent  maybe  less than the one produced by  the  market noises, that is, $\mathcal G_t \subseteq \mathcal F_t$, where $\mathbb G\equiv\{\mathcal G_t\}_{t\geq0}$ is  the  information available to the agent. There are considerable literatures on related topics \cite{Hu2008Partial,Zhen2010A,huang2009a,Wang2015A,Wang2017An,wang2018a}.  In particular, Huang et al.  \cite{huang2009a} derived a   necessary condition  for optimality  of BSDE with partial information and applied their results to two classes of  LQ problems.   Wang et al. \cite{Wang2015A} and  Wang  et al. \cite{Wang2017An} concerned LQ problems with partially observable information driven by  FBSDE and mean field FBSDE, respectively. LQ non-zero sum stochastic differential game of BSDE is considered in Wang et al.  \cite{wang2018a}.      They   obtained   feedback Nash equilibrium points by FBSDE and Riccati equation under  asymmetric information.


Our work distinguishes itself from existing literatures in the following aspects. (i) Both the generator of dynamic system and the cost functional contain diffusion terms $Z_1$ and $Z_2$.  Moreover, our results are obtained under some usual conditions (see Assumptions $A1$ and $A2$ in section 2). In the literatures on this topic, diffusion terms $Z_1$ and $Z_2$ are usually assumed not to be contained into the generator (see \cite{huang2009a}, \cite{Wang2017An}), or there are some additional conditions to ensure the solvability of Riccati equation (see \cite{Wang2015A}, \cite{wang2018a}).
(ii) Sufficient and necessary conditions of optimality are established, which provide an expression  for optimal control via the solution of   stochastic Hamiltonian system.  (iii)  Explicit representations for optimal control   in terms of  three Riccati equations, a  BSDE with filtering, and an SDE with filtering are obtained, as well as  the associated optimal  cost. The derivation of associated Riccati equations is extremely different from \cite{Lim2001Linear} and \cite{Li2017Linear}, since the stochastic Hamiltonian system is an FBSDE with filtering.  Moreover,   the uniqueness and existence of solution to   BSDE \eqref{BSDE partial} is first obtained, which is important in deriving  explicit  representations  for optimal control and associated optimal cost. (iv) Last but not least, we consider two special scalar-valued control problems of BSDEs with partial  information. In the case of $H=N_1=0$, we obtain explicit solutions of the stochastic Hamiltonian system, as well as related Riccati equations. In the case of $C_2=0$, we give some numerical simulations to illustrate our theoretical results.

The rest of this paper is organized as follows. In Section 2,
we formulate the stochastic  backward LQ control problem and give some preliminary results.  Section 3 aims to   decouple the associated stochastic Hamiltonian system and derive some Riccati equations. In Section 4, we give    explicit  representations of optimal control and the associated optimal cost. Section 5 is devoted to solving two special scalar-valued control problems and giving some numerical simulations. Finally, we conclude this paper.

\section{Preliminaries}
Let $(\Omega,\mathcal F, \mathbb F, \mathbb P)$ be a complete filtered probability space and let $T>0$ be a fixed time horizon. Let $\{(W_{1t}, W_{2t}): 0\leq t\leq T\}$ be  a $\mathbb R^2$-valued standard Wiener process, defined on $(\Omega,\mathcal F, \mathbb F, \mathbb P)$. $\mathbb F\equiv\{\mathcal F_t\}_{t\geq0}$ is a natural filtration of $(W_{1}, W_{2})$ augmented by all   $\mathbb P$-null sets. Let $\mathcal F_t^{\beta}=\sigma\{\beta_{s},\ 0\leq s\leq t\}$ be  the filtration generated by a stochastic process $\beta$. Let $\mathbb R^{n\times m}$ be the set of all $n\times m$   matrices and $\mathbb S^n$ be the set  of all $n\times n$ symmetric matrices. For a matrix $M\in \mathbb R^n$, let $M^\top$ be its transpose. The inner product $\langle\cdot,\cdot\rangle$ on $\mathbb R^{n\times m}$ is defined by $\langle M,N\rangle\mapsto tr(M^\top N)$ with an  induced norm  $|M|=\sqrt{tr(M^\top M)}$. In particular, we denote  by $\mathbb S_+^n$ ($\widehat{\mathbb S}_+^n$) the set  of all $n\times n$ (uniformly) positive definite matrices.  For any Euclidean space $M$, we adopt the following notations:\\
$\mathcal L_{\mathcal F_T}^2(\Omega; M)=\Big\{\zeta: \Omega\to M| \zeta $ is an $\mathcal F_T$-measurable random variable, $\mathbb E[|\zeta|^2]<\infty$\Big\};\\
$\mathcal L^\infty(0,T;M)=\Big\{v:[0,T]\to M|v$ is a bounded function\Big\};\\
$\mathcal L_{\mathcal F}^2(0,T;M)=\Big\{v:[0,T]\times \Omega\to M|v$  is an $\{\mathcal F_t\}_{t\geq0}$-adapted stochastic process, $\mathbb E\left[\int_0^T|v_t|^2dt\right]<\infty\Big\};$\\
$\mathcal S_{\mathcal F}^2(0,T;M)=\Big\{v:[0,T]\times \Omega\to  M|v$ is an $\{\mathcal F_t\}_{t\geq0}$-adapted stochastic process and has continuous paths, $\mathbb E\left[\sup_{t\in[0,T]}|v_t|^2\right]<\infty\Big\}$.

Consider a controlled linear BSDE
\begin{equation}\label{station 1}
\left\{ \begin{aligned}
dY_t=&\big(A_tY_t+B_tv_t+C_{1t}Z_{1t}+C_{2t}Z_{2t}\big)dt+Z_{1t}dW_{1t}+Z_{2t}dW_{2t}, \quad  t\in{[0,T]},\\
Y_T=&\ \zeta,
\end{aligned}\right.
\end{equation}
where $\zeta\in L_{\mathcal F_T}^2(\Omega; \mathbb R^n)$ and $v$, valued in $\mathbb R^m$, is a control process.
Introduce an  admissible control set\\
$\mathcal V_{ad}[0,T]=\Big\{v:[0,T]\times \Omega\to\mathbb R^m|v\ is\  \{\mathcal F_t^{W_1}\}_{t\geq0}-$ adapted, $\mathbb E\left[\int_0^T|v_t|^2dt\right]<\infty\Big\}.$\\
Any $v\in \mathcal V_{ad}[0,T]$ is called an  admissible control.\\
\textbf{Assumption $A1$:} The coefficients of  dynamic system satisfy
\begin{equation*}
A, C_1,  C_2\in\mathcal L^\infty(0,T;\mathbb R^{n\times n}),\
 B \in\mathcal L^\infty(0,T;\mathbb R^{n\times m}).
\end{equation*}
Under Assumption $A1$, dynamic system   \eqref{station 1} admits a unique solution pair $(Y, Z_1, Z_2)\in\mathcal S_{\mathcal F}^2(0,T;\mathbb R^n)\times\mathcal L_{\mathcal F}^2(0,T;\mathbb R^n)\times\mathcal L_{\mathcal F}^2(0,T;\mathbb R^n)$, which is called the corresponding state process, for any $v\in \mathcal V_{ad}[0,T]$ (see Pardoux and Peng \cite{Pardoux1990Adapted}, Yong and Zhou \cite{Yong1999}).
We introduce a quadratic cost functional
\begin{equation}\label{cost functional}
\begin{aligned}
J(v)=&\frac{1}{2}\mathbb E\Bigg[Y_0^\top GY_0+\int_0^T\Big(Y_t^\top H_tY_t+v_t^\top R_tv_t+Z_{1t}^\top N_{1t}Z_{1t}+Z_{2t}^\top N_{2t}Z_{2t}\Big)dt\Bigg].
\end{aligned}
\end{equation}
\textbf{Assumption $A2$:} The weighting matrices in  cost functional satisfy
\begin{equation*}
\begin{aligned}
&H, N_1,  N_2\in\mathcal L^\infty(0,T;\mathbb S_+^{n}),\ R\in\mathcal L^\infty(0,T;\widehat{\mathbb S}_+^m), G\in \mathbb S_+^{n}.
\end{aligned}
\end{equation*}
Our stochastic backward  LQ control problem   can be stated as follows.
\\\textbf{Problem BLQ}. Find a $v^*\in \mathcal V_{ad}[0,T]$ such that
\begin{equation}\label{inf}
J(v^*)=\inf_{v\in \mathcal V_{ad}[0,T]}J(v).
\end{equation}
Any $v^*\in \mathcal V_{ad}[0,T]$ satisfying \eqref{inf} is called an optimal control, and the state process $(Y^*,Z_1^*,Z_2^*)$   is called an optimal state process. Under Assumptions $A1$ and $A2$, Problem BLQ is   uniquely solvable for any terminal state $\zeta\in\mathcal L_{\mathcal F_T}^2(\Omega; \mathbb R^n)$ (see Li et al. \cite{Li2017Linear}).  We  suppressed the time argument  in the sequel of this paper wherever necessary, for the sake of notation simplicity. The following theorem is a necessary condition of optimality, which is easy to be   obtained from Theorem 3.1 in Huang et al. \cite{huang2009a}.
\begin{theorem}\label{necessary condition}
Under  Assumptions $A1-A2$, if $v^*$ is an optimal control of Problem BLQ and $(Y^*, Z_1^*, Z_2^*)$ is  the corresponding optimal state process, then
\begin{equation}\label{adjoint equation}
\left\{\begin{aligned}dX^*=&-\left(A^\top X^*+HY^*\right)dt-\left(C_{1}^\top X^*+N_{1}Z_{1}^*\right)dW_{1}-\left(C_{2}^\top X^*+N_{2}Z_{2}^*\right)dW_{2}, \\X_{0}^*=&-GY_{0}^*
\end{aligned}
\right.
\end{equation}
admits a unique solution such that
$$\mathbb E\left[R_tv_t^*+B_t^\top X_t^*|\mathcal F^{W_1}\right]=0,\ t\in [0,T],\ a.s..$$
\end{theorem}
According to the above analysis, we end up with a  stochastic Hamiltonian system
\begin{equation}\label{H   system}
\left\{\begin{aligned}
&dY=\left(AY+Bv+C_{1}Z_{1}+C_{2}Z_{2}\right)dt+Z_{1}dW_{1}+Z_{2}dW_{2}, \\&dX=-\left(A^\top X+HY\right)dt-\left(C_{1}^\top X+N_{1}Z_{1}\right)dW_1
-\left(C_{2}^\top X+N_{2}Z_{2}\right)dW_{2},\\&Y=\zeta,\ \ \ \  X_{0}=-GY_{0},
\\& \mathbb E[R_tv_t+B_t^\top X_t|\mathcal F_t^{W_1}]=0.
\end{aligned}
\right.
\end{equation}
This is a coupled FBSDE with filtering. Note that the coupling comes from the last equation in \eqref{H   system}, which is also called a stationarity condition. We point out that in our setting, the stationarity condition involves a conditional expectation, which makes the decoupling of this stochastic Hamiltonian system different and difficult. We now prove the sufficiency of the above result.
\begin{theorem}\label{sufficient condition}
Let Assumption $A1 - A2$ hold. If $(X^*, Y^*, Z_1^*, Z_2^*,  v^*)$ is an adapted solution to  stochastic Hamiltonian system  \eqref{H   system}, then $v^*$ is an optimal control.
\end{theorem}
\begin{proof}
For any $v\in \mathcal V_{ad}[0,T]$, let $({Y},{Z_1},{Z_2})$ be the corresponding state process. Let $(\widetilde{Y},\widetilde{Z_1},\widetilde{Z_2})$ satisfies
\begin{equation*}
\left\{\begin{aligned}
d\widetilde{Y}=&\left[A\widetilde{Y}+B(v-v^*)+C_{1}\widetilde{Z_{1}} +C_{2}\widetilde{Z_{2}}\right]dt+\widetilde{Z_{1}}dW_{1}+\widetilde{Z_{2}}dW_{2}, \\\widetilde{Y}_T=&\ 0.
\end{aligned}
\right.
\end{equation*}
According to the existence and uniqueness of solution to   BSDE, we have  $\widetilde{Y}=Y-Y^*, \widetilde{Z_1}=Z_1-Z_1^*, \widetilde{Z_2}=Z_2-Z_2^*$.
With the notation, we derive
\begin{align*}
J(v)-J(v^*)=\mathbb E\bigg[ Y_0^{*\top} G\widetilde Y_0
+\int_0^T\Big({Y^*}^\top H\widetilde Y+{v^*}^\top R(v-v^*)+{Z_{1}^*}^\top N_{1}\widetilde Z_{1}+{Z_{2}^*}^\top N_{2}\widetilde Z_{2}\Big)dt\bigg]+\widetilde{J},
\end{align*}
where
\begin{align*}
\widetilde{J}=&\frac{1}{2}\mathbb E\bigg[\widetilde{Y_0}^\top G\widetilde Y_0+\int_0^T\Big(\widetilde{Y}^\top H\widetilde Y+ (v-v^*)^\top R (v-v^*)+\widetilde{Z_{1}}^\top N_{1}\widetilde Z_{1}+\widetilde{Z_{2}}^\top N_{2}\widetilde Z_{2}\Big)dt\bigg].
\end{align*}
It is easy to see that $\widetilde{J}\geq 0$ under Assumption $A2$.
Further,
\begin{align*}
\ \mathbb E\left[Y_0^{*\top} G\widetilde Y_0\right]
=&\ \mathbb E\bigg[\int_0^T\Big(\langle A\widetilde{Y}+B(v-v^*)+C_{1}\widetilde{Z_{1}}+C_{2}\widetilde{Z_{2}},X^*\rangle-\langle \widetilde{Y}, A^\top X^*+HY^*\rangle\\&-\langle \widetilde{Z_{1}}, C_{1}^\top X^*+N_{1}Z_{1}^*\rangle-\langle \widetilde{Z_{2}}, C_{2}^\top X^*+N_{2}Z_{2}^*\rangle \Big)dt\bigg]\\=&\ \mathbb E\bigg[\int_0^T\Big(\langle v-v^*,B^\top X^*\rangle-\langle \widetilde{Y},HY^*\rangle-\langle \widetilde{Z_{1}}, N_{1}Z_{1}^*\rangle-\langle \widetilde{Z_{2}}, N_{2}Z_{2}^*\rangle \Big)dt\bigg].
\end{align*}
Thus, we have
\begin{align*}
J(v)-J(v^*)=\mathbb E\left[\int_0^T\langle v-v^*,Rv^*+B^\top X^*\rangle dt\right]+\widetilde{J}\geq 0.
\end{align*}
Then, $v^*$ is an optimal control.
\end{proof}

\section{Decoupling  stochastic Hamiltonian system \eqref{H   system}}
 In this section, we use the decoupling method for general FBSDE introduced in \cite{Ma2007Forward} to solve    stochastic Hamiltonian system \eqref{H   system}, which is an FBSDE with filtering. Different from the results in \cite{Lim2001Linear},   we  obtain three Riccati equations, an   BSDE with filtering and an  SDE with filtering.  For simplicity of notation, we denote $\widehat \beta_{t}=\mathbb E[\beta_{t}|\mathcal F_{t}^{W_1}]$.
To  be precise, we assume   that
\begin{equation}\label{duality partial 1}
Y=\Upsilon\widehat X+\varphi,
\end{equation}
where $\Upsilon$ is a differential and deterministic matrix-valued function with a  terminal condition $\Upsilon_T=0$, and $\varphi$ is a stochastic process satisfying the BSDE
\begin{equation*}
\left\{\begin{aligned}
d\varphi=&\ \lambda dt +\eta_{1}dW_{1}+\eta_{2}dW_{2},
\\\varphi_T=&\ \zeta,
\end{aligned}
\right.
\end{equation*}
for   $\{\mathcal F_t\}_{t\geq0}$-adapted processes $\lambda$, $\eta_1$ and $\eta_2$. According to Theorem 2.1 in Wang et al. \cite{Wangbook} (see also Theorem 5.7 in Xiong \cite{Xiong2008An}  and Theorem 8.1 in  Liptser and Shiryayev \cite{Liptser} ), we have
\begin{equation*}
\left\{\begin{aligned}
d\widehat X=&-\left(A^\top\widehat X+H\widehat Y\right)dt-\left(C_{1}^\top\widehat X+N_{1}\widehat Z_{1}\right)dW_{1},
\\\widehat X_0=&-G\widehat Y_0.
\end{aligned}
\right.
\end{equation*}
Applying It\^o formula to \eqref{duality partial 1}, we  get
\begin{equation*}
\begin{aligned}
0=&\ dY-\dot\Upsilon\widehat Xdt-\Upsilon d\widehat X-d\varphi\\=&\left(AY+Bv+C_{1}Z_{1}+C_{2}Z_{2}\right)dt +Z_{1}dW_{1}+Z_{2}dW_{2}-\dot \Upsilon\widehat Xdt+\Upsilon\left(A^\top\widehat X+H\widehat Y\right)dt\\&+\Upsilon(C_{1}^\top\widehat X+N_{1}\widehat Z_{1})dW_{1}-\lambda dt -\eta_{1}dW_{1}-\eta_{2}dW_{2}.
\end{aligned}
\end{equation*}
This implies
\begin{equation} \label{relationback}
\left\{
\begin{aligned}
&AY-BR^{-1}B^\top\widehat X+C_{1}Z_{1}+C_{2}Z_{2}-\dot \Upsilon\widehat X+\Upsilon\left(A^\top\widehat X+H\widehat Y\right)-\lambda=0,\\
&Z_{1}+\Upsilon(C_{1}^\top\widehat X+N_{1}\widehat {Z_{1}})-\eta_{1}=0,\\
&Z_{2}-\eta_{2}=0.
\end{aligned}
\right.
\end{equation}
Assuming that $I+\Upsilon N_1$ is invertible, we  have
\begin{equation}\label{z in terms of eta}
\left\{
\begin{aligned}
&Z_{1}=\eta_{1}-\widehat\eta_{1}+(I+\Upsilon N_{1})^{-1}(\widehat\eta_{1}-\Upsilon C_{1}^\top\widehat X),\\
&Z_{2}=\eta_{2}.
\end{aligned}
\right.
\end{equation}
Substituting \eqref{duality partial 1} and \eqref{z in terms of eta} into the first equation in \eqref{relationback}, we obtain
\begin{align*}
&\ A(\Upsilon\widehat X+\varphi)-BR^{-1}B^\top\widehat X+C_{1}(\eta_{1}-\widehat\eta_{1})+C_{1}(I+\Upsilon N_{1})^{-1}(\widehat\eta_{1}-\Upsilon C_{1}^\top\widehat X)\\&+C_{2}\eta_{2}-\dot\Upsilon\widehat X+\Upsilon A^\top\widehat X+\Upsilon H(\Upsilon\widehat X+\widehat \varphi)-\lambda=0.
\end{align*}
Then $\Upsilon$ satisfies a Riccati equation
\begin{equation}\label{Riccati Sigma partial}
\left\{
\begin{aligned}
&\dot\Upsilon-\Upsilon A^\top-A\Upsilon-\Upsilon H\Upsilon+BR^{-1}B^\top+C_{1}(I+\Upsilon N_{1})^{-1}\Upsilon C_{1}^\top=0,\\
&\Upsilon_T=0,
\end{aligned}
\right.
\end{equation}
 and $\varphi$ satisfies a  BSDE
\begin{equation}\label{BSDE partial}
\left\{
\begin{aligned}
d\varphi=&\Big[A\varphi+\Upsilon H\widehat \varphi+C_{1}(\eta_{1}-\widehat\eta_{1})+C_{1}(I+\Upsilon N_{1})^{-1}\widehat\eta_{1}+C_{2}\eta_{2}\Big]dt\\&+\eta_{1}dW_{1}+\eta_{2}dW_{2},\\
\varphi_T=&\ \zeta.
\end{aligned}
\right.
\end{equation}
Riccati equation \eqref{Riccati Sigma partial} admits a unique solution $\Upsilon \in \mathcal L^\infty(0,T;\mathbb S_+^{n})$ under Assumptions $A1-A2$ (see \cite{Li2017Linear,Lim2001Linear}). Note that \eqref{BSDE partial} is a BSDE with filtering, for which the solvability  has not been given in literatures before. We will specified this problem in Section 4. In order to give the optimal control with  a feedback representation, we conjecture  that
\begin{equation}\label{duality partial 2}
X=-\Gamma_{1}(Y-\widehat Y)-\Gamma_{2}\widehat Y-\psi,
\end{equation}
where $\Gamma_1$ and $\Gamma_2$ are differential and deterministic matrix-valued functions with  initial conditions  $\Gamma_{10}=G$ and  $\Gamma_{20}=G$, respectively;  $\psi$ is a stochastic process satisfying an SDE
\begin{equation*}
\left\{\begin{aligned}
d\psi =&\ \alpha_{0} dt +\alpha_{1} dW_{1} +\alpha_{2} dW_{2} ,
\\\psi_0=&\ 0,
\end{aligned}
\right.
\end{equation*}
where   $\alpha_0 $, $\alpha_1$ and $\alpha_2$ are   $\{\mathcal F_t\}_{t\geq 0}$-adapted processes. Note that
\begin{equation*}
\left\{\begin{aligned}
d\widehat Y=&\left(A\widehat Y-BR^{-1}B^\top\widehat X+C_{1}\widehat Z_{1}+C_{2}\widehat Z_{2}\right)dt+\widehat Z_{1}dW_{1},
\\\widehat Y_T=&\ \widehat\zeta,
\end{aligned}
\right.
\end{equation*}
where $\widehat \zeta=\mathbb E[\zeta|\mathcal F_T^{W_1}]$.
Hence,
\begin{equation*}
\left\{\begin{aligned}
d(Y-\widehat Y)=&\left[A(Y-\widehat Y)+C_{1}(Z_{1}-\widehat Z_{1})+C_{2}(Z_{2}-\widehat Z_{2})\right]dt+(Z_{1}-\widehat Z_{1})dW_{1}+Z_{2}dW_{2},
\\Y_T-\widehat Y_T=&\ \zeta-\widehat\zeta.
\end{aligned}
\right.
\end{equation*}
Applying It\^o formula to \eqref{duality partial 2}, we obtain
\begin{equation*}
\begin{aligned}
0=&\ dX+\dot \Gamma_{1}(Y-\widehat Y)dt+\Gamma_{1}d(Y-\widehat Y)+\dot \Gamma_{2}\widehat Ydt+\Gamma_{2}d\widehat Y+d\psi\\
=&-\left(A^\top X+HY\right)dt-\left(C_{1}^\top X+N_{1}Z_{1}\right)dW_{1}
-\left(C_{2}^\top X+N_{2}Z_{2}\right)dW_{2}\\&+\dot \Gamma_{1}(Y-\widehat Y)dt+\Gamma_{1}\Big[A(Y-\widehat Y)+C_{1}(Z_{1}-\widehat Z_{1})+C_{2}(Z_{2}-\widehat Z_{2})\Big]dt+\Gamma_{1}(Z_{1}-\widehat Z_{1})dW_1\\&+\Gamma_{1}Z_{2}dW_{2}+\dot \Gamma_{2}\widehat Ydt+\Gamma_{2}\Big(A\widehat Y-BR^{-1}B^\top\widehat X+C_{1}\widehat Z_{1}+C_{2}\widehat Z_{2}\Big)dt+\Gamma_{2}\widehat Z_{1}dW_{1}\\&+\alpha_{0}dt +\alpha_{1}dW_{1}+\alpha_{2}dW_{2}.
\end{aligned}
\end{equation*}
It yields
\begin{equation*}
\left\{
\begin{aligned}
&-\left(A^\top X+HY\right)+\dot \Gamma_{1}(Y-\widehat Y)+\Gamma_{1}\Big[A(Y-\widehat Y)+C_{1}(Z_{1}-\widehat Z_{1})+C_{2}(Z_{2}-\widehat Z_{2})\Big]\\&+\dot \Gamma_{2}\widehat Y+\Gamma_{2}\Big(A\widehat Y-BR^{-1}B^\top\widehat X+C_{1}\widehat Z_{1}+C_{2}\widehat Z_{2}\Big)+\alpha_{0}=0,\\
&-(C_{1}^\top X+N_{1}Z_{1})+\Gamma_{1}(Z_{1}-\widehat Z_{1})+\Gamma_{2}\widehat Z_{1}+\alpha_{1}=0,\\
&-(C_{2}^\top X+N_{2}Z_{2})+\Gamma_{1}Z_{2}+\alpha_{2}=0.
\end{aligned}
\right.
\end{equation*}
Assuming that $I+\Gamma_{2}\Upsilon$ is invertible,  we arrive at
\begin{equation*}
\left\{
\begin{aligned}
\alpha_{1}=&\ (N_{1}-\Gamma_{1})(\eta_{1}-\widehat \eta_{1})+(N_{1}-\Gamma_{2})(I+\Upsilon N_{1})^{-1}\left[\widehat \eta_{1}+\Upsilon C_{1}^\top(I+\Gamma_{2}\Upsilon)^{-1}(\Gamma_{2}\widehat \varphi+\widehat\psi)\right]\\&-C_{1}^\top \Gamma_{1}(\varphi-\widehat \varphi)-C_{1}^\top(\psi-\widehat \psi)-C_{1}^\top(I+\Gamma_{2}\Upsilon)^{-1}(\Gamma_{2}\widehat \varphi+\widehat\psi),
\\\alpha_{2}=&\ (N_{2}-\Gamma_{1})\eta_{2}-C_{2}^\top \Gamma_{1}(\varphi-\widehat \varphi)-C_{2}^\top(\psi-\widehat \psi)-C_{2}^\top(I+\Gamma_{2}\Upsilon)^{-1}(\Gamma_{2}\widehat \varphi+\widehat\psi).
\end{aligned}
\right.
\end{equation*}
Further, it follows from \eqref{z in terms of eta} and \eqref{duality partial 2} that
\begin{align*}
&A^\top \Gamma_{1}(Y-\widehat Y)+A^\top\Gamma_{2}\widehat Y+A^\top\psi-HY+\dot \Gamma_{1}(Y-\widehat Y)+\Gamma_{1}A(Y-\widehat Y)\\&+\Gamma_{1}C_{1}(\eta_{1}-\widehat\eta_{1})+\Gamma_{1}C_{2}(\eta_{2}-\widehat\eta_{2})+\dot \Gamma_{2}\widehat Y+\Gamma_{2}A\widehat Y+\Gamma_{2}BR^{-1}B^\top(\Gamma_{2}\widehat Y+\widehat\psi)\\&+\Gamma_{2}C_{1}(I+\Upsilon N_{1})^{-1}\left[\widehat\eta_{1}+\Upsilon C_{1}^\top(\Gamma_{2}\widehat Y+\widehat\psi)\right]+\Gamma_{2}C_{2}\widehat\eta_{2}+\alpha_{0}=0.
\end{align*}
Introduce
\begin{equation}\label{Riccati M1 partial}
\left\{
\begin{aligned}
&\dot \Gamma_{1}+\Gamma_{1}A+A^\top \Gamma_{1}-H=0,\\
&\Gamma_{10}=G,
\end{aligned}
\right.
\end{equation}
\begin{equation}\label{Riccati M2 partial}
\left\{
\begin{aligned}
&\dot \Gamma_{2}+\Gamma_{2}A+A^\top \Gamma_{2}+\Gamma_{2}BR^{-1}B^\top \Gamma_{2}+\Gamma_{2}C_{1}(I+\Upsilon N_{1})^{-1}\Upsilon C_{1}^\top \Gamma_{2}-H=0,\\
&\Gamma_{20}=G,
\end{aligned}
\right.
\end{equation}
and
\begin{equation}\label{SDE partial}
\left\{
\begin{aligned}
d\psi=&-\Big[A^\top\psi+\Gamma_{2}BR^{-1}B^\top\widehat\psi+\Gamma_{2}C_{1}(I+\Upsilon N_{1})^{-1}\left(\widehat \eta_{1}+\Upsilon C_{1}^\top\widehat\psi\right)\\&+\Gamma_{2}C_{2}\widehat\eta_{2}+\Gamma_{1}C_{1}(\eta_{1}-\widehat\eta_{1})+\Gamma_{1}C_{2}(\eta_{2}-\widehat\eta_{2})\Big]dt
\\&+\Big\{(N_{1}-\Gamma_{1})(\eta_{1}-\widehat \eta_{1})-C_{1}^\top \Gamma_{1}(\varphi-\widehat \varphi)-C_{1}^\top(\psi-\widehat \psi)\\&+(N_{1}-\Gamma_{2})(I+\Upsilon N_{1})^{-1}\left[\widehat \eta_{1}+\Upsilon C_{1}^\top(I+\Gamma_{2}\Upsilon)^{-1}(\Gamma_{2}\widehat \varphi+\widehat\psi)\right]\\&-C_{1}^\top(I+\Gamma_{2}\Upsilon)^{-1}(\Gamma_{2}\widehat \varphi+\widehat\psi)\Big\}dW_{1}\\&+\Big\{(N_{2}-\Gamma_{1})\eta_{2}-C_{2}^\top \left[\Gamma_{1}(\varphi-\widehat \varphi)+(\psi-\widehat \psi)\right]-C_{2}^\top(I+\Gamma_{2}\Upsilon)^{-1}(\Gamma_{2}\widehat \varphi+\widehat\psi)\Big\}dW_{2},\\
\psi_0=&\ 0.
\end{aligned}
\right.
\end{equation}
There is   a unique solution $\Gamma_{1}\in\mathcal  L^\infty(0,T;\mathbb S_+^{n})$ to Riccati equation \eqref{Riccati M1 partial}, since  Assumptions $A1$ and $A2$ hold (see Yong and Zhou \cite{Yong1999}). Corollary 4.6 in Lim and Zhou \cite{Lim2001Linear} implies that  \eqref{Riccati M2 partial} admits a unique solution $\Gamma_{2}\in\mathcal  L^\infty(0,T;\mathbb S_+^{n})$.  Once $\Upsilon$, $\Gamma_{1}$, $\Gamma_{2}$ and the solution  $(\varphi, \eta_{1}, \eta_{2})$ of \eqref{BSDE partial} are known, the solvability of  \eqref{SDE partial} will be obtained immediately.
\section{Explicit representations of optimal control and optimal cost}
Now we  would like to   give   explicit formulas  of optimal control and associated optimal cost  in terms of  Riccati equations \eqref{Riccati Sigma partial}, \eqref{Riccati M1 partial}, \eqref{Riccati M2 partial},   BSDE \eqref{BSDE partial} and   SDE \eqref{SDE partial}. We first    prove that \eqref{BSDE partial}    admits a unique solution. Consider a BSDE
\begin{equation}\label{general BSDE partial}
\left\{
\begin{aligned}
dP=&\ g(t,P,Q_{1}, Q_{2}, \widehat{P},\widehat{Q_{1}}, \widehat{Q_{2}})dt+Q_{1}dW_{1}+Q_{2}dW_{2},\\P_T=&\ \zeta,
\end{aligned}
\right.
\end{equation}
where $\widehat P_{t}=\mathbb E[P_{t}|\mathcal F_{t}^{W_{1}}], \widehat Q_{1t}=\mathbb E[Q_{1t}|\mathcal F_{t}^{W_{1}}], \widehat Q_{2t}=\mathbb E[Q_{2t}|\mathcal F_{t}^{W_{1}}]$.\\
We assume that\\
\textbf{Assumption $A3$:} There exists a constant $L$, such that, $\mathbb P$-a.s., for all $t\in [0,T]$, $p,q_{1}, q_{2},\bar{p},\bar{q_{1}}, \bar{q_{2}} $, $p',q_{1}', q_{2}',\bar{p}',\bar{q_{1}}', \bar{q_{2}}'\in \mathbb R^n$,
\begin{align*}
&\left|g(t,p,q_{1}, q_{2},\bar{p},\bar{q_{1}}, \bar{q_{2}})-g(t,p',q_{1}', q_{2}',\bar{p}',\bar{q_{1}}', \bar{q_{2}}')\right|\\\leq
&\ L\big(|p-p'|+|q_{1}-q_{1}'|+|q_{2}-q_{2}'|+|\bar{p}-\bar{p}'|+|\bar{q_{1}}-\bar{q_{1}}'|+|\bar{q_{2}}-\bar{q_{2}}'|\big).\end{align*}
\textbf{Assumption $A4$:} $g(\cdot,0,0,0,0,0,0)\in\mathcal L_{\mathcal F}^2(0,T;\mathbb R^n)$.\\
\begin{lemma}
Let  Assumptions $A3$  and $A4$ hold. For any $\zeta\in\mathcal L_{\mathcal F}^2(\Omega;\mathbb R^n)$,    BSDE \eqref{general BSDE partial} admits a unique solution $(P,Q_{1}, Q_{2})\in \mathcal S_{\mathcal F}^2(0,T;\mathbb R^n)\times\mathcal L_{\mathcal F}^2(0,T;\mathbb R^n)\times\mathcal L_{\mathcal F}^2(0,T;\mathbb R^n)$.
\end{lemma}
\begin{proof}
We first introduce a norm on $\mathcal L_{\mathcal F}^2(0,T;\mathbb R^{n+n+n})$, which is equivalent to the canonical norm
$$||u||_{\delta}=\left(\mathbb E\left[\int_0^T|u_t|^2e^{\delta t}dt\right]\right)^{\frac{1}{2}},\ \delta>0.$$
The parameter $\delta$ will be specified later.
For any $(p,q_{1}, q_{2})\in L_{\mathcal F}^2(0,T;\mathbb R^{n+n+n})$, the following BSDE
\begin{equation*}
\left\{
\begin{aligned}
dP=&\ g(t,P,Q_{1}, Q_{2}, \widehat{p},\widehat{q_{1}}, \widehat{q_{2}})dt+Q_{1}dW_{1}+Q_{2}dW_{2}, \\P_T=&\ \zeta
\end{aligned}
\right.
\end{equation*}
admits a unique solution $(P,Q_{1}, Q_{2})\in\mathcal L_{\mathcal F}^2(0,T;\mathbb R^{n+n+n})$.
We then introduce a mapping $(P,Q_{1}, Q_{2})=\mathbf I(p,q_{1}, q_{2})$: $\mathcal L_{\mathcal F}^2(0,T;\mathbb R^{n+n+n}) \to\mathcal L_{\mathcal F}^2(0,T;\mathbb R^{n+n+n})$ by
\begin{equation*}
\left\{
\begin{aligned}
dP=&\ g(t,P,Q_{1}, Q_{2}, \widehat{p},\widehat{q_{1}}, \widehat{q_{2}})dt+Q_{1}dW_{1}+Q_{2}dW_{2}, \\P_T=&\ \zeta.
\end{aligned}
\right.
\end{equation*}
For any $(p,\ q_{1},\ q_{2})$, $(p',\ q_{1}',\ q_{2}')\in L_{\mathcal F}^2(0,T;\mathbb R^{n+n+n})$, we denote $(P,Q_{1}, Q_{2})=\mathbf I(p,q_{1}, q_{2})$, $(P',Q_{1}', Q_{2}')=\mathbf I(p',q_{1}', q_{2}')$, $(\widetilde{p},\ \widetilde{q_{1}},\ \widetilde{q_{2}})=(p-p', q_{1}-q_{1}', q_{2}-q_{2}')$  and $(\widetilde{P},\ \widetilde{Q_{1}},\ \widetilde{Q_{2}})=(P-P', Z_{1}-Z_{1}', Q_{2}-Q_{2}')$. Applying It\^o formula to $|\widetilde{P}_t|^2e^{\delta t}$  and taking conditional expectations, we get
\begin{align*}
&\ |\widetilde{P}_t|^2+\mathbb E\left[\int_t^T\delta e^{\delta (s-t)}|\widetilde{P_s}|^2ds\Big|\mathcal F_t\right]+\mathbb E\left[\int_t^T e^{\delta (s-t)}(|\widetilde{Q_{1s}}|^2+|\widetilde{Q_{2s}}|^2)ds\Big|\mathcal F_t\right]
\\=&\ 2\mathbb E\bigg[\int_t^Te^{\delta (s-t)}\langle\widetilde{P}_{s},g(s,P',Q_{1}', Q_{2}', \widehat{p'},\widehat{q_{1}'}, \widehat{q_{2}'})-g(s,P,Q_{1}, Q_{2}, \widehat{p},\widehat{q_{1}}, \widehat{q_{2}})\rangle ds\Big|\mathcal F_t\bigg]
\\\leq&\ 2L\mathbb E\bigg[\int_t^Te^{\delta (s-t)}|\widetilde{P}_{s}|\Big(|\widetilde{P}_{s}|+|\widetilde{Q}_{1 s}|+|\widetilde{Q}_{2s}|+|\widehat{\widetilde{p}}_{s}| +|\widehat{\widetilde{q}}_{1 s}|+|\widehat{\widetilde{q}}_{2s}|\Big)ds\Big|\mathcal F_t\bigg]
\\\leq&\ \mathbb E\bigg[\int_t^Te^{\delta (s-t)}\Big((2L+4L^2+\frac{\delta}{2})|\widetilde{P}_{s}|^2+\frac{1}{2}|\widetilde{Q}_{1 s}|^2+\frac{1}{2}|\widetilde{Q}_{2s}|^2\Big)ds\Big|\mathcal F\bigg]\\&+\mathbb E\bigg[\int_t^Te^{\delta (s-t)}\frac{6L^2}{\delta}\Big(|\widehat{\widetilde{p}}_{s}|^2+|\widehat{\widetilde{q}}_{1 s}|^2+|\widehat{\widetilde{q}}_{2s}|^2\Big)ds\Big|\mathcal F_t\bigg].
\end{align*}
Thus we have
\begin{align*}
&\left(\frac{\delta}{2}-2L-4L^2\right)\mathbb E\left[\int_0^Te^{\delta s}|\widetilde{P}_{s}|^2ds\right]+\frac{1}{2}\mathbb E\left[\int_0^Te^{\delta s}\Big(|\widetilde{Q_{1}}_{s}|^2+|\widetilde{Q_{2}}_{s}|^2\Big)ds\right]\\\leq&\ \frac{6L^2}{\delta}\mathbb E\left[\int_0^Te^{\delta s}\Big(|\widehat{\widetilde{p}}_{s}|^2+|\widehat{\widetilde{q_{1}}}_{s}|^2+|\widehat{\widetilde{q_{2}}}_{s}|^2\Big)ds\right].
\end{align*}
Taking $\delta=24L^2+4L+1$, we arrive at
\begin{align*}
&\mathbb E\left[\int_0^Te^{\delta s}\Big(|\widetilde{P_{s}}|^2+|\widetilde{Q_{1s}}|^2+|\widetilde{Q_{2s}}|^2\Big)ds\right]\\\leq &\ \frac{1}{2}\mathbb E\left[\int_0^Te^{\delta s}\Big(|\widetilde{p_{s}}|^2+|\widetilde{q}_{1s}|^2+|\widetilde{q}_{2s}|^2\Big)ds\right],
\end{align*}
which implies  $||(\widetilde{P}, \widetilde{Q_{1}}, \widetilde{Q_{2}})||_{\delta}\leq \frac{1}{\sqrt{2}}||(\widetilde{p}, \widetilde{q_{1}}, \widetilde{q_{2}})||_{\delta}$. That is, $\mathbf I$ is a contraction on $\mathcal L_{\mathcal F}^2(0,T;\mathbb R^{n+n+n})$, endowed with the norm $||\cdot||_{\delta}$. According to the contraction mapping theorem, we know that there is a unique fixed point $(P,Q_{1}, Q_{2})\in\mathcal L_{\mathcal F}^2(0,T;\mathbb R^{n+n+n})$, such that $\mathbf I(P,Q_{1}, Q_{2})=(P,Q_{1}, Q_{2})$, which is exactly the solution of \eqref{general BSDE partial}. We now proceed to prove that $P \in\mathcal S_{\mathcal F}^2(0,T;\mathbb R^{n})$.   Using 	
Jensen inequality, H\"{o}lder's inequality and  Burkholder-Davis-Gundy's inequality yields
\begin{align*}
&\ \mathbb E\left[\sup_{t\in [0,T]}\Big|P_t\Big|^2\right]\\\leq&\ 4\mathbb E\left[|\zeta|^2\right]+4\mathbb E\left[\sup_{t\in [0,T]}\Big|\int_t^Tg(s,P,Q_{1}, Q_{2}, \widehat{Q},\widehat{Q_{1}}, \widehat{Q_{2}})ds\Big|^2\right]\\&+4\mathbb E\left[\sup_{t\in [0,T]}\Big|\int_t^TQ_{1s}dW_{1s}\Big|^2\right]+4\mathbb E\left[\sup_{t\in [0,T]}\Big|\int_t^TQ_{2s}dW_{2s}\Big|^2\right]\\\leq&\ 4\mathbb E\left[|\zeta|^2\right]+4T\mathbb E\left[\sup_{t\in [0,T]}\left(\int_t^T|g(s,P,Q_{1}, Q_{2}, \widehat{P},\widehat{Q_{1}}, \widehat{Q_{2}})|^2ds\right)\right]\\&+16\mathbb E\left[\int_0^T\big|Q_{1t}\big|^2 dt\right]+16\mathbb E\left[\int_0^T\big|Q_{2t}\big|^2 dt\right]\\<&\ \infty.
\end{align*}
Therefore, we obtain $P\in\mathcal  S_{\mathcal F}^2(0,T;\mathbb R^{n})$.
\end{proof}
\begin{remark}
Equation \eqref{BSDE partial} is a linear BSDE with filtering, where   the   generator  satisfies Assumptions  $A3-A4$. Then it follows   that \eqref{BSDE partial} admits a unique solution $(\varphi,\eta_{1},\eta_{2})\in\mathcal S_{\mathcal F}^2(0,T;\mathbb R^n)\times \mathcal L_{\mathcal F}^2(0,T;\mathbb R^n)\times\mathcal L_{\mathcal F}^2(0,T;\mathbb R^n)$.
\end{remark}

We have the following theorem which specifies the solvability of stochastic Hamiltonian system \eqref{H system} and gives some relations between the forward component and the backward components.
\begin{theorem}\label{Relation}
Under  Assumptions  $A1-A2$,   stochastic Hamiltonian system \eqref{H system} admits  a unique solution $(X,Y,Z_{1},Z_{2},v)$. Moreover, we have the following relations
\begin{equation*}
\left\{\begin{aligned}
&Y=\Upsilon\widehat X+\varphi,\\
&Z_{1}=\eta_{1}-\widehat\eta_{1}+(I+\Upsilon N_{1})^{-1}(\widehat\eta_{1}-\Upsilon C_{1}^\top\widehat X),\\
&Z_{2}=\eta_{2},
\\&v=-R^{-1}B^\top\widehat X,
\\&Y_0=(I+\Upsilon_0G)^{-1}\varphi_0,
\end{aligned}
\right.
\end{equation*}
where $\Upsilon$ and $(\varphi,\eta_{1},\eta_{2})$ are  solutions to  \eqref{Riccati Sigma partial} and \eqref{BSDE partial}, respectively.
\end{theorem}
\begin{proof}
Consider the following SDE with filtering
\begin{equation}\label{SDE for p}
\left\{\begin{aligned}
d\bar X=&-\left[A^\top\bar X+H(\Upsilon\widehat {\bar X}+\varphi)\right]dt\\&-\Big[C_{1}^\top\bar X+N_{1}(\eta_{1}-\widehat\eta_{1})+N_{1}(I+\Upsilon N_{1})^{-1}\big(\widehat\eta_{1}-\Upsilon C_{1}^\top\widehat {\bar X}\big)\Big]dW_{1}
\\&-\left(C_{2}^\top\bar X+N_{2}\eta_{2}\right)dW_{2},\\\bar X_0=&-(I+G\Upsilon_0)^{-1}G\varphi_0,
\end{aligned}
\right.
\end{equation}
where $\Upsilon$ and $(\varphi,\eta_{1},\eta_{2})$ are  solutions of \eqref{Riccati Sigma partial} and \eqref{BSDE partial}, respectively.
According to Theorem 2.1 in Wang et al. \cite{Wangbook}, we get
\begin{equation}\label{SDE for filtering p}
\left\{\begin{aligned}
d\widehat{\bar X}=&-\left(A^\top\widehat{\bar X}+H\Upsilon \widehat {\bar X}+H\widehat{\varphi}\right)dt \\&-\Big[(I+N_{1}\Upsilon)^{-1}C_{1}^\top\widehat{\bar X}+N_{1}(I+\Upsilon N_{1})^{-1}\widehat\eta_{1}\Big]dW_{1},
\\\widehat{\bar X}_0=&-(I+G\Upsilon_0)^{-1}G\varphi_0.
\end{aligned}
\right.
\end{equation}
From the theory of linear SDE, \eqref{SDE for filtering p} has a unique solution $\widehat {\bar X}$. Then it follows that \eqref{SDE for p} also admits a unique solution $\bar X$. We define
$$\bar Y=\Upsilon\widehat {\bar X}+\varphi.$$
By using It\^o formula, $\bar Y$ satisfies
\begin{equation*}
\begin{aligned}
d\bar Y=&\Big[\Upsilon A^\top+A\Upsilon+\Upsilon H\Upsilon-BR^{-1}B^\top-C_{1}(I+\Upsilon N_{1})^{-1}\Upsilon C_{1}^\top\Big]\widehat {\bar X}dt\\&-\Upsilon\left(A^\top\widehat{\bar X}+H\Upsilon \widehat {\bar X}+H\widehat{\varphi}\right)dt -\Upsilon\Big[(I+N_{1}\Upsilon)^{-1}C_{1}^\top\widehat{\bar X}+N_{1}(I+\Upsilon N_{1})^{-1}\widehat\eta_{1}\Big]dW_{1}\\&+\Big[A\varphi+\Upsilon H\widehat {\varphi}+C_{1}(\eta_{1}-\widehat\eta_{1})+C_{1}(I+\Upsilon N_{1})^{-1}\widehat\eta_{1}+C_{2}\eta_{2}\Big]dt+\eta_{1}dW_{1}+\eta_{2}dW_{2}
\\=&\Big[A\bar Y-BR^{-1}B^\top\widehat{\bar X}+C_{1}(\eta_{1}-\widehat\eta_{1})+C_{1}(I+\Upsilon N_{1})^{-1}(\widehat\eta_{1}-\Upsilon C_{1}^\top\widehat {\bar X})+C_{2}\eta_{2}\Big]dt\\&+\left[\eta_{1}-\widehat\eta_{1}+(I+\Upsilon N_{1})^{-1}(\widehat\eta_{1}-\Upsilon C_{1}^\top\widehat{\bar X})\right]dW_{1}+\eta_{2}dW_{2},
\end{aligned}
\end{equation*}
with an initial condition $Y_0=(I+\Upsilon_0G)^{-1}\varphi_0$. Defining 

 \begin{equation*}
\left\{\begin{aligned}
&\bar Z_{1}=\eta_{1}-\widehat\eta_{1}+(I+\Upsilon N_{1})^{-1}(\widehat\eta_{1}-\Upsilon C_{1}^\top\widehat{\bar X}),\\
&\bar Z_{2}=\eta_{2},\\
&\bar v=-R^{-1}B^\top\widehat {\bar X}.
\end{aligned}
\right.
\end{equation*}
It is obvious  that $(\bar X,\bar Y,\bar Z_{1},\bar Z_{2}, \bar v)$ is a solution to stochastic Hamiltonian system \eqref{H system}.

We now turn to prove the uniqueness. Suppose that equation \eqref{H system} admits two solutions  $(X,Y,Z_{1},Z_{2},v)$ and $(X',Y',Z_{1}',Z_{2}',v')$, respectively. Let $(\widetilde{X},\widetilde{Y},\widetilde{Z_{1}},\widetilde{Z_{2}},\widetilde{v})
=(X-X',Y-Y',Z_{1}-Z_{1}',Z_{2}-Z_{2}',v-v')$. Thus $(\widetilde{X},\widetilde{Y},\widetilde{Z_{1}},\widetilde{Z_{2}},\widetilde{v})$ satisfies
\begin{equation*}
\left\{\begin{aligned}
&d\widetilde{Y}=\left(A\widetilde{Y}+B\widetilde{v}+C_{1}\widetilde{Z_{1}}+C_{2}\widetilde{Z_{2}}\right)dt +\widetilde{Z_{1}}dW_{1}+\widetilde{Z_{2}}dW_{2}, \\&d\widetilde{X}=-\left(A^\top\widetilde{X}+H\widetilde{Y}\right)dt-\left(C_{1}^\top\widetilde{X}+N_{1}\widetilde{Z_{1}}\right)dW_{1}
-\left(C_{2}^\top\widetilde{X}+N_{2}\widetilde{Z_{2}}\right)dW_{2},\\&\widetilde{Y}_T=0,\ \ \ \  \widetilde{X}_0=-G\widetilde{Y}_0,\\&\mathbb E[R_t\widetilde{v}_t+B_t^\top\widetilde{X}_t|\mathcal F_t^{W_{1}}]=0.
\end{aligned}
\right.
\end{equation*}
Applying It\^o formula to $\widetilde{Y}^\top\widetilde{X}$, we obtain
\begin{equation*}
\begin{aligned}
&\mathbb E\left[\widetilde{Y}_0^\top G\widetilde{Y}_0\right]\\=&-\mathbb E\left[\int_0^T\widetilde{Y}^\top\left(A^\top\widetilde{X}+H\widetilde{Y}\right)dt\right]+\mathbb E\left[\int_0^T\left(A\widetilde{Y}-BR^{-1}B^\top\widehat {\widetilde{X}}+C_{1}\widetilde{Z_{1}}+C_{2}\widetilde{Z_{2}}\right)^\top\widetilde{X}dt\right]
\\&-\mathbb E\left[\int_0^T\widetilde{Z_{1}}^\top\left(C_{1}^\top\widetilde{X}+N_{1}\widetilde{Z_{1}}\right)dt\right]
-\mathbb E\left[\int_0^T\widetilde{Z_{2}}^\top\left(C_{2}^\top\widetilde{X}+N_{2}\widetilde{Z_{2}}\right)dt\right]
\\=&-\mathbb E\Big[\int_0^T\big(\widetilde{Y}^\top H\widetilde{Y}+\widehat {\widetilde{X}}^\top BR^{-1}B^\top \widetilde{X}+\widetilde{Z_{1}}^\top N_{1}\widetilde{Z_{1}}+\widetilde{Z_{2}}^\top N_{2}\widetilde{Z_{2}}\big)dt\Big].
\end{aligned}
\end{equation*}
We adopt the same procedure as in the proof of Theorem \ref{sufficient condition}. Since $G, H, R, N_{1}, N_{2}$ satisfy Assumption $A2$, it follows that
\begin{equation*}
\begin{aligned}
\mathbb E\left[\int_0^T\widehat {\widetilde{X}}^\top BR^{-1}B^\top \widetilde{X}dt\right]=0.
\end{aligned}
\end{equation*}
Recalling $R$ is uniformly positive, it yields $$B_t^\top \widehat{\widetilde{X}_t}=0,\ \  \ a.e.\ t\in[0,T],\ \  \mathbb P-a.s..$$ With the equality, $(\widetilde{Y},\widetilde{Z_{1}},\widetilde{Z_{2}})$ satisfies
\begin{equation}\label{unique y}
\left\{\begin{aligned}
&d\widetilde{Y}=\left(A\widetilde{Y}+C_{1}\widetilde{Z_{1}}+C_{2}\widetilde{Z_{2}}\right)dt +\widetilde{Z_{1}}dW_{1}+\widetilde{Z_{2}}dW_{2}, \\&\widetilde{Y}_T=0.
\end{aligned}
\right.
\end{equation}
It is easy to see that \eqref{unique y} admits a unique solution $(\widetilde{Y},\widetilde{Z_{1}},\widetilde{Z_{2}})\equiv0$. Then
\begin{equation*}
\left\{\begin{aligned}
&d\widetilde{X}=-A^\top\widetilde{X}dt-C_{1}^\top\widetilde{X}dW_{1}
-C_{2}^\top\widetilde{X}dW_{2},\\&\widetilde{X}_0=-G\widetilde{Y}_0.
\end{aligned}
\right.
\end{equation*}
Hence it follows from the uniqueness of solution that  $\widetilde{X}\equiv0$. The proof is completed.
\end{proof}%

To summarize the above  analysis, we establish the following main result.
\begin{theorem}\label{main partial}
Let Assumptions $A1-A2$ hold and let $\zeta\in\mathcal  L_{\mathcal F}^2(\Omega;\mathbb R^n)$ be given. Let $\Upsilon$, $\Gamma_{1}$, $\Gamma_{2}$ be the   solutions of Riccati equations \eqref{Riccati Sigma partial}, \eqref{Riccati M1 partial} and \eqref{Riccati M2 partial}, respectively.   Let $(\varphi,\eta_{1},\eta_{2})$ and $\psi$ be the   solutions of \eqref{BSDE partial} and \eqref{SDE partial}, respectively. Then the   BSDE with filtering
\begin{equation*}
\left\{\begin{aligned}
dY=&\ \Big(AY+BR^{-1}B^\top \Gamma_{2}\widehat {Y}+BR^{-1}B^\top\widehat \psi+C_{1}Z_{1}+C_{2}Z_{2}\Big)dt+Z_{1}dW_{1}+Z_{2}dW_{2},\\Y_T=&\ \zeta
\end{aligned}
\right.
\end{equation*}
admits a unique solution  $(Y,Z_{1},Z_{2})$. By defining
\begin{equation*}
\left\{\begin{aligned}
X=&-\Gamma_{1}(Y-\widehat Y)-\Gamma_{2}\widehat Y-\psi,\\
v=&\ R^{-1}B^\top \Gamma_{2}\widehat {Y}+R^{-1}B^\top\widehat \psi,
\end{aligned}
\right.
\end{equation*}
the 5-tuple $(X,Y,Z_{1},Z_{2},v)$ is an adapted  solution to FBSDE \eqref{H system} and $v$ is an optimal control of Problem BLQ. The corresponding optimal cost is
\begin{equation}\label{cost partial}
\begin{aligned}
J(v)=&\frac{1}{2}\mathbb E\left[\langle\widehat\zeta, \Sigma_T\widehat\zeta\rangle\right]+\frac{1}{2}\mathbb E\left[\int_0^T\Big(\langle H\varphi,\varphi\rangle- \langle\widehat{\varphi},H \widehat{\varphi}\rangle\Big) dt\right]\\&+\frac{1}{2}\mathbb E\left[\int_0^T\big\langle \big[N_{1}(I+\Upsilon N_{1})^{-1}-\Sigma\big]\widehat \eta_{1},\widehat \eta_{1}\big\rangle\right]\\&+\frac{1}{2}\mathbb E\left[\int_0^T\Big(\langle N_{1}(\eta_{1}-\widehat \eta_{1}),\eta_{1}-\widehat \eta_{1}\rangle +\langle N_{2}\eta_{2},\eta_{2} \rangle\Big) dt\right]\\&-\mathbb E\left[\int_0^T\Big\langle\widehat \varphi, \Sigma\Big(C_{1}(I+\Upsilon N_{1})^{-1}\widehat{\eta_{1}}+C_{2}\widehat\eta_{2}\Big)\Big\rangle dt\right],
\end{aligned}
\end{equation}
where $\Sigma$ is the solution of
\begin{equation*}\label{cost partial P}
\left\{\begin{aligned}
&\dot \Sigma+\Sigma(A+\Upsilon H)+(A+\Upsilon H)^\top \Sigma-H=0,\\&\Sigma_0=G(I+\Upsilon_0G)^{-1}.
\end{aligned}
\right.
\end{equation*}
\end{theorem}

\begin{proof}
We need only to prove \eqref{cost partial}.
Substituting \eqref{duality partial 1}, \eqref{z in terms of eta} into the cost functional, we derive
\begin{align*}
J(v)=&\ \frac{1}{2}\big\langle G(I+\Upsilon_0G)^{-1}\varphi_0, (I+\Upsilon_0G)^{-1}\varphi_0\big\rangle+\frac{1}{2}\mathbb E\left[\int_0^T\Big\langle\Upsilon\widehat X+\varphi, H(\Upsilon\widehat X+\varphi)\Big\rangle dt\right]
\\&+\frac{1}{2}\mathbb E\left[\int_0^T\Big\langle R^{-1}B^\top\widehat{X},B^\top\widehat{X}\Big\rangle dt\right]+\frac{1}{2}\mathbb E\left[\int_0^T\langle\eta_{2},N_{2}\eta_{2} \rangle dt\right]
\\&+\frac{1}{2}\mathbb E\bigg[\int_0^T\Big\langle \eta_{1}-\widehat\eta_{1}+(I+\Upsilon N_{1})^{-1}(\widehat\eta_{1}-\Upsilon C_{1}^\top\widehat X),N_{1}\Big[\eta_{1}-\widehat\eta_{1}+(I+\Upsilon N_{1})^{-1}(\widehat\eta_{1}-\Upsilon C_{1}^\top\widehat X)\Big]\Big\rangle dt\bigg]
\\=&\ \frac{1}{2}\big\langle G(I+\Upsilon_0G)^{-1}\varphi_0, (I+\Upsilon_0G)^{-1}\varphi_0\big\rangle
\\&+\frac{1}{2}\mathbb E\bigg[\int_0^T\Big\langle\widehat X, \Big[\Upsilon H\Upsilon+BR^{-1}B^\top+C_{1}(I+\Upsilon N_{1})^{-1}\Upsilon N_{1}\Upsilon(I+N_{1}\Upsilon)^{-1}C_{1}^\top\Big]\widehat X\Big\rangle dt\bigg]
\\&+\mathbb E\left[\int_0^T\Big\langle \Upsilon H \varphi,\widehat X\Big\rangle dt\right]-\mathbb E\bigg[\int_0^T\Big\langle\eta_{1}-\widehat \eta_{1}+(I+\Upsilon N_{1})^{-1}\widehat \eta_{1}, N_{1}(I+\Upsilon N_{1})^{-1}\Upsilon C_{1}^\top\widehat X\Big\rangle dt\bigg]\\&+\frac{1}{2}\mathbb E\bigg[\int_0^T\Big\langle\eta_{1}-\widehat \eta_{1}+(I+\Upsilon N_{1})^{-1}\widehat \eta_{1}, N_{1}\Big[\eta_{1}-\widehat \eta_{1}+(I+\Upsilon N_{1})^{-1}\widehat \eta_{1}\Big]\Big\rangle dt\bigg]\\&+\frac{1}{2}\mathbb E\left[\int_0^T\Big(\langle\eta_{2}, N_{2}\eta_{2} \rangle+\langle \varphi,H\varphi\rangle\Big) dt\right].
\end{align*}

Applying  It\^o formula to $\langle\widehat X, \Upsilon \widehat X\rangle$, we have
\begin{align*}
&\big\langle(I+G\Upsilon_0)^{-1}G\varphi_0, \Upsilon_0(I+G\Upsilon_0)^{-1}G\varphi_0\big\rangle
\\=&-\mathbb E\left[\int_0^T\Big\langle\widehat X, \Upsilon\big[(A+\Upsilon H)^\top\widehat X+Q\widehat\varphi\big]\Big\rangle dt\right]-\mathbb E\left[\int_0^T\Big\langle(A+\Upsilon H)^\top\widehat X+H\widehat \varphi, \Upsilon\widehat X\Big\rangle dt\right]
\\&+\mathbb E\bigg[\int_0^T\Big\langle\widehat X, \Big(\Upsilon A^\top+A\Upsilon+\Upsilon H\Upsilon-BR^{-1}B^\top-C_{1}(I+\Upsilon N_{1})^{-1}\Upsilon C_{1}^\top\Big)\widehat X\Big\rangle dt\bigg]
\\&+\mathbb E\bigg[\int_0^T\Big\langle (I+N_{1}\Upsilon)^{-1}C_{1}^\top\widehat{X}+N_{1}(I+\Upsilon N_{1})^{-1}\widehat\eta_{1},\Upsilon\Big[(I+N_{1}\Upsilon)^{-1}C_{1}^\top\widehat{X}+N_{1}(I+\Upsilon N_{1})^{-1}\widehat\eta_{1}\Big]\Big\rangle dt\bigg]
\\=&-\mathbb E\bigg[\int_0^T\Big\langle\widehat X, \Big[\Upsilon H\Upsilon+BR^{-1}B^\top+C_{1}(I+\Upsilon N_{1})^{-1}\Upsilon N_{1}\Upsilon(I+N_{1}\Upsilon)^{-1}C_{1}^\top\Big]\widehat X\Big\rangle dt\bigg]
\\&+\mathbb E\left[\int_0^T\Big\langle N_{1}(I+\Upsilon N_{1})^{-1}\widehat{\eta_{1}},\Upsilon N_{1}(I+\Upsilon N_{1})^{-1}\widehat{\eta_{1}}\Big\rangle dt\right]\\&+2\mathbb E\left[\int_0^T\Big\langle (I+N_{1}\Upsilon)^{-1}C_{1}^\top\widehat X, \Upsilon  N_{1}(I+\Upsilon N_{1})^{-1}\widehat{\eta_{1}}\Big\rangle dt\right]-2\mathbb E\left[\int_0^T\Big\langle \Upsilon H\widehat \varphi,\widehat X\Big\rangle dt\right].
\end{align*}
With the equality, we derive
\begin{equation*}
\begin{aligned}
J(v)=&\ \frac{1}{2}\langle (I+G\Upsilon_0)^{-1}G\varphi_0, \varphi_0\rangle+\frac{1}{2}\mathbb E\left[\int_0^T\langle H\varphi,\varphi\rangle dt\right]\\&+\frac{1}{2}\mathbb E\bigg[\int_0^T\Big(\langle N_{1}(I+\Upsilon N_{1})^{-1}\widehat \eta_{1},\widehat \eta_{1}\rangle+ \langle N_{1}(\eta_{1}-\widehat \eta_{1}),\eta_{1}-\widehat \eta_{1}\rangle+\langle N_{2}\eta_{2},\eta_{2} \rangle \Big) dt\bigg].
\end{aligned}
\end{equation*}
Recalling  that $ \varphi$ satisfies \eqref{BSDE partial} and  applying It\^o formula to $\langle\widehat \varphi, \Sigma\widehat \varphi\rangle$, we have
\begin{align*}
\langle (I+G\Upsilon_0)^{-1}G\varphi_0, \varphi_0\rangle
=&\ \mathbb E[\langle\widehat\zeta, \Sigma_T\widehat\zeta\rangle]
-\mathbb E\left[\int_0^T\Big(\langle\widehat{\varphi},H \widehat{\varphi}\rangle+\langle\widehat{\eta_{1}},\Sigma \widehat{\eta_{1}}\rangle\Big) dt\right]\\&-2\mathbb E\left[\int_0^T\Big\langle\widehat\varphi, \Sigma\Big[C_{1}(I+\Upsilon N_{1})^{-1}\widehat{\eta_{1}}+C_{2}\widehat\eta_{2}\Big]\Big\rangle dt\right].
\end{align*}
We obtain
\begin{align*}
J(v)=&\frac{1}{2}\mathbb E\left[\langle\widehat\zeta, \Sigma_T\widehat\zeta\rangle\right]+\frac{1}{2}\mathbb E\left[\int_0^T\Big(\langle H\varphi,\varphi\rangle- \langle\widehat{\varphi}, H \widehat{\varphi}\rangle\Big) dt\right]\\&+\frac{1}{2}\mathbb E\left[\int_0^T\big\langle \big[N_{1}(I+\Upsilon N_{1})^{-1}-\Sigma\big]\widehat \eta_{1},\widehat \eta_{1}\big\rangle\right]\\&+\frac{1}{2}\mathbb E\left[\int_0^T\Big(\langle N_{1}(\eta_{1}-\widehat \eta_{1}),\eta_{1}-\widehat \eta_{1}\rangle +\langle N_{2}\eta_{2},\eta_{2} \rangle\Big) dt\right]\\&-\mathbb E\left[\int_0^T\big\langle\widehat \varphi, \Sigma\big[C_{1}(I+\Upsilon N_{1})^{-1}\widehat{\eta_{1}}+C_{2}\widehat\eta_{2}\big]\big\rangle dt\right].
\end{align*}
Then our claims follow.
\end{proof}
\begin{remark}
When we consider the complete information case, i. e., $W_{2}$ disappears in \eqref{station 1}. Let  $\zeta$ be  an $\mathcal F_T^{W_{1}}$-measurable square integrable   random variable. Let  $v$ be  an $\{\mathcal F_t^{W_{1}}\}_{t\geq0}-$adapted and  square integrable stochastic process. Then Theorem \ref{main partial} coincides with Theorem 3.2 in Lim and Zhou \cite{Lim2001Linear}.
\end{remark}
\begin{remark}
In Huang et al. \cite{2020A}, an optimal control for Problem BLQ with  feedback representation is given. We point out that their results rely   on the condition that the solution of \eqref{BSDE partial} satisfies $\eta_{2}=0$.
\end{remark}
\section{One-dimensional case}
In this section, we consider  two scalar-valued backward LQ problems with partial information and give more detailed analyses. In the case of   $H=N_{1}=0$, we work out an  explicit control problem and show the detailed procedure to obtain the feedback representation of   optimal control using our theoretical results. In the case of   $C_{2}=0$, we give some numerical simulations to illustrate our theoretical results, since we can not obtain   explicit solutions of related stochastic Hamiltonian system and Riccati equation.

\subsection{Special case: $H=N_{1}=0$}
Under Assumptions $A1$ and $A2$, let all the coefficients of \eqref{station 1} and \eqref{cost functional} are constants, and $$\zeta=e^{(a-\frac{1}{2}b^2-\frac{1}{2}c^2)T+bW_{1}+cW_{2}}.$$ In this case, \eqref{station 1} is given by
\begin{equation*}\label{station App1}
\left\{ \begin{aligned}
dY_t=&\left(AY_t+Bv_t+C_{1}Z_{1t}+C_{2}Z_{2t}\right)dt+Z_{1t}dW_{1t}+Z_{2t}dW_{2t}, \qquad t\in{[0,T]},\\
Y_T=&\ \zeta.
\end{aligned}\right.
\end{equation*}
The cost functional takes the  form of
\begin{equation*}\label{cost functional App1}
\begin{aligned}
J(v)=&\ \frac{1}{2}\mathbb E\left[GY_0^2+\int_0^T(Rv^2+N_{2}Z_{2}^2)dt\right].
\end{aligned}
\end{equation*}
Then Problem BLQ is   stated as follows.\\
\textbf{Problem BLQA}. Find a $v^* \in \mathcal V[0,T]$ such that
\begin{equation*}
J(v^*)=\inf_{v \in \mathcal V[0,T]}J(v),
\end{equation*}
where the admissible control set is given by\\
$\mathcal V[0,T]=\Big\{v:[0,T]\times \Omega\to \mathbb R|v \ is\  \{\mathcal F_t^{W_{1}}\}_{t\geq0}$-adapted, $\mathbb E\left[\int_0^Tv_t^2dt\right]<\infty\Big\}.$\\

The corresponding stochastic Hamiltonian system reads
\begin{equation}\label{H   system App1}
\left\{\begin{aligned}
&dY_t=\left(AY_t+Bv_t+C_{1}Z_{1t}+C_{2}Z_{2t}\right)dt+Z_{1t}dW_{1t}+Z_{2t}dW_{2t}, \\&dX_t=-AX_tdt-C_{1}X_tdW_{1t}
-(C_{2}X_t+N_{2}Z_{2t})dW_{2t},\\&Y_T=\zeta,\ \ \ \  X_0=-GY_0,
\\& \mathbb E[Rv_t+BX_t|\mathcal F_t^{W_{1}}]=0.
\end{aligned}
\right.
\end{equation}
We introduce
\begin{equation*}\label{Riccati Sigma partial APP1}
\left\{
\begin{aligned}
&\dot\Upsilon_t-(2A-C_{1}^2)\Upsilon_t+\frac{B^2}{R}=0,\\
&\Upsilon_T=0,
\end{aligned}
\right.
\end{equation*}
and
\begin{equation*}\label{BSDE partial App1}
\left\{
\begin{aligned}
d\varphi_t=&\big(A\varphi_t+C_{1}\eta_{1t}+C_{2}\eta_{2t}\big)dt+\eta_{1t}dW_{1t}+\eta_{2t}dW_{2t},\\
\varphi_T=&\ \zeta.
\end{aligned}
\right.
\end{equation*}
It is easy to see that
\begin{equation*}\label{solu App1 sigma}
\Upsilon_t=\left\{
\begin{array}{ll}
\frac{B^2}{R(2A-C_{1}^2)}\Big(1-e^{(2A-C_{1}^2)(t-T)}\Big),  &2A-C_{1}^2\neq 0,\\
\frac{B^2(T-t)}{R}, &2A-C_{1}^2=0,
\end{array}
\right.
\end{equation*}
and
\begin{equation}\label{solu App1 BSDE}
\left\{
\begin{aligned}
\varphi_t=&\exp\Big[(a-bC_{1}-cC_{2}-A)T\\&+(A+bC_{1}+cC_{2}-\frac{1}{2}b^2-\frac{1}{2}c^2)t+bW_{1t}+cW_{2t}\Big],\\
\eta_{1t}=&\ b\varphi_t,\\
\eta_{2t}=&\ c\varphi_t.
\end{aligned}
\right.
\end{equation}
Taking $t=0$ in \eqref{solu App1 BSDE},  we have
\begin{equation*}
\varphi_0=\exp\Big[(a-bC_{1}-cC_{2}-A)T\Big].
\end{equation*}
Then it follows from \eqref{H   system App1} and \eqref{solu App1 BSDE} that
\begin{equation*}
\left\{\begin{aligned}
dX_t=&-AX_tdt-C_{1}X_tdW_{1t}
-(C_{2}X_t+N_{2}\eta_{2t})dW_{2t},\\X_0=&-(I+G\Upsilon_0)^{-1}G\varphi_0,
\end{aligned}
\right.
\end{equation*}
which admits a unique solution
\begin{equation}\label{solu adjoint}
\begin{aligned}
X_t=&\Psi_t\Big[X_0-\int_0^t\Psi_{s}^{-1}C_{2}N_{2}\eta_{2s}ds-\int_0^t\Psi_{s}^{-1} N_{2}\eta_{2s}dW_{2s}\Big],
\end{aligned}
\end{equation}
with
\begin{equation*}
\Psi_t=\exp\left[-\big(A+\frac{1}{2}C_{1}^2+\frac{1}{2}C_{2}^2\big)t-C_{1}W_{1t}-C_{2}W_{2t}\right].
\end{equation*}
Further,
\begin{equation}\label{control App1nondeedback}
\begin{aligned}
v_t=-R^{-1}B\widehat X_t=-R^{-1}BX_0\exp\Big[-\big(A+\frac{1}{2}C_{1}^2\big)t-C_{1}W_{1t}\Big].
\end{aligned}
\end{equation}
Theorem \ref{sufficient condition} implies that $v$ given by \eqref{control App1nondeedback}  is an optimal control of Problem BLQA.

In the following, we aim to derive a  feedback representation of $v$. For this end, we introduce
\begin{equation}\label{Riccati M1 partial App1}
\left\{
\begin{aligned}
&\dot \Gamma_{1t}+2A\Gamma_{1t}=0,\\
&\Gamma_{10}=G,
\end{aligned}
\right.
\end{equation}
\begin{equation}\label{Riccati M2 partial App1}
\left\{
\begin{aligned}
&\dot \Gamma_{2t}+2A\Gamma_{2t}+\left(\frac{B^2}{R}+C_{1}^2\Upsilon_t\right)\Gamma_{2t}^2=0,\\
&\Gamma_{20}=G,
\end{aligned}
\right.
\end{equation}
and
\begin{equation*}
\left\{
\begin{aligned}
d\psi_t=&-\Big[A\psi_t+\Big(\frac{B^2\Gamma_{2t}}{R}+C_{1}^2\Gamma_{2t}\Upsilon_t\Big)\widehat\psi_t +C_{1}\Gamma_{2t}\widehat \eta_{1t}\\&+C_{1}\Gamma_{1t}(\eta_{1t}-\widehat\eta_{1t})+C_{2}\Gamma_{1t}(\eta_{2t}-\widehat\eta_{2t}) +C_{2}\Gamma_{2t}\widehat\eta_{2t}\Big]dt\\&
-\Big[\Gamma_{1t}(\eta_{1t}-\widehat \eta_{1t})+\Gamma_{2t}\widehat \eta_{1t}+C_{1}(\Gamma_{2t}\widehat \varphi_t+\widehat\psi_t)\\&+C_{1}\Gamma_{1t}(\varphi_t-\widehat \varphi_t)+C_{1}(\psi_t-\widehat \psi_t)\Big]dW_{1t}\\&-\Big[(\Gamma_{1t}-N_{2})\eta_{2t}+C_{2}\Gamma_{1t}(\varphi_t-\widehat \varphi_t)\\&+C_{2}(1+\Gamma_{2t}\Upsilon_t)^{-1}(\Gamma_{2t}\widehat \varphi_t+\widehat\psi_t)+C_{2}(\psi_t-\widehat \psi_t)\Big]dW_{2t},\\
\psi_{0}=&\ 0.
\end{aligned}
\right.
\end{equation*}
Solving  \eqref{Riccati M1 partial App1} and \eqref{Riccati M2 partial App1}, we get
\begin{equation*}\label{solu App1 M1}
\Gamma_{1t}=Ge^{-2At},
\end{equation*}
and
\begin{equation*}\label{solu App1 M2}
\Gamma_{2t}=\frac{G\exp(-2At)}{1+G\int_0^t\exp(-2As)(\frac{B^2}{R}+C_{1}^2\Upsilon_{s})ds},
\end{equation*}
respectively.\\
According to Theorem 2.1 in Wang et al. \cite{Wangbook}, we have
\begin{equation*}\label{SDE partial App1}
\left\{
\begin{aligned}
d\widehat\psi_t=&-\Big(\mathcal A_t\widehat\psi_t+\mathcal B_t\Big)dt
-\Big(C_{1}\widehat\psi_t+\mathcal D_t\Big)dW_{1t},\\
\widehat\psi_0=&\ 0,
\end{aligned}
\right.
\end{equation*}
where
\begin{equation*}
\left\{
\begin{aligned}
\mathcal A_t=&A+\frac{B^2\Gamma_{2t}}{R}+C_{1}^2\Gamma_{2t}\Upsilon_t,
\\\mathcal B_t=&C_{1}\Gamma_{2t}\widehat \eta_{1t}+C_{2}\Gamma_{2t}\widehat\eta_{2t},
\\\mathcal D_t=&\Gamma_{2t}\widehat \eta_{1t}+C_{1}\Gamma_{2t}\widehat \varphi_t.
\end{aligned}
\right.
\end{equation*}
Similarly, we derive
\begin{equation*}
\begin{aligned}
\widehat\psi_t=&\Phi_t\bigg[-\int_0^t\Phi^{-1}_{s}\Big(\mathcal B_{s}+C_{1}\mathcal D_{s}\Big)ds-\int_0^t\Phi^{-1}_{s}\mathcal D_{s}dW_{1s}\bigg],
\end{aligned}
\end{equation*}
where
\begin{equation*}
\Phi_t=\exp\left[-\int_0^t\big(\mathcal A_{s}+\frac{1}{2}C_{1}^2\big)ds-C_{1}W_{1t}\right].
\end{equation*}
Then Theorem \ref{main partial} implies that  \eqref{control App1nondeedback} admits a feedback representation below
\begin{equation*}
v_t=R^{-1}B\Gamma_{2t}\widehat {Y}_t+R^{-1}B\widehat \psi_t,
\end{equation*}
where $Y$ satisfies
\begin{equation}\label{state APP1}
\left\{\begin{aligned}
dY_t=&\ \Big(AY_t+B^2R^{-1}\Gamma_{2t}\widehat {Y}_t+B^2R^{-1}\widehat \psi_t+C_{1}Z_{1t}+C_{2}Z_{2t}\Big)dt+Z_{1t}dW_{1t}+Z_{2t}dW_{2t},\\Y_T=&\ \zeta.
\end{aligned}
\right.
\end{equation}
The corresponding optimal cost  is
\begin{equation*}
\begin{aligned}
J(v)=&\frac{G\varphi_0^2}{2+2G\Upsilon_0}+\frac{1}{2}\mathbb E\left[\int_0^T\langle N_{2}\eta_{2t},\eta_{2t} \rangle dt\right].
\end{aligned}
\end{equation*}
\begin{remark}
It follows from Theorem \ref{Relation}  that the solution of \eqref{state APP1} is given by
\begin{equation*}
\left\{\begin{aligned}
&Y_t=\Upsilon_t\widehat X_t+\varphi_t,\\
&Z_{1t}=\eta_{1t}- C_{1}\Upsilon_t\widehat X_t,\\
&Z_{2t}=\eta_{2t}.
\end{aligned}
\right.
\end{equation*}
Here, $X$ and $(\varphi,\eta_{1},\eta_{2})$ are given by \eqref{solu adjoint} and \eqref{solu App1 BSDE}, respectively.  Note that equation \eqref{state APP1} is a BSDE with filtering, which is difficult to obtain the explicit solution in general.
\end{remark}
\subsection{Special case: $C_{2}=0$}
In this case, \eqref{station 1} is written as
\begin{equation*}\label{station App2}
\left\{ \begin{aligned}
dY=&\left(AY+Bv+C_{1}Z_{1}\right)dt+Z_{1}dW_{1}+Z_{2}dW_{2}, t\in{[0,T]},\\
Y_T=&\ \zeta.
\end{aligned}\right.
\end{equation*}
Cost functional \eqref{cost functional} takes the   form of
\begin{equation*}\label{cost functional App2}
\begin{aligned}
J(v)=&\frac{1}{2}\mathbb E\Bigg[GY_0^2+\int_0^T\Big(HY^2+Rv^2+N_{1} Z_{1}^2+N_{2} Z_{2}^2 \Big)dt\Bigg].
\end{aligned}
\end{equation*}
Then Problem BLQ is formulated as follows.\\
\textbf{Problem BLQB}. Find a $v^*\in \mathcal V[0,T]$ such that
\begin{equation*}
J(v^*)=\inf_{v\in \mathcal V[0,T]}J(v).
\end{equation*}

The corresponding stochastic Hamiltonian system reads
\begin{equation*}\label{H   system App2}
\left\{\begin{aligned}
&dY=\left(AY+Bv+C_{1}Z_{1}\right)dt+Z_{1}dW_{1}+Z_{2}dW_{2}, \\&dX=-\left(AX+HY\right)dt-\left(C_{1}X+N_{1}Z_{1}\right)dW_{1}
-N_{2}Z_{2}dW_{2},\\&Y=\zeta,\ \ \ \  X_0=-GY_0,
\\& \mathbb E[R_tv_t+B_tX_t|\mathcal F_t^{W_{1}}]=0.
\end{aligned}
\right.
\end{equation*}
According to  Theorem \ref{main partial}, the optimal control is
$$v_t=R_t^{-1}B_t\Gamma_{2t}\widehat {Y}_t+R_t^{-1}B_t\widehat \psi_t$$  with
\begin{equation*}
\left\{\begin{aligned}
dY=&\left[AY+B(R^{-1}B\Gamma_{2}\widehat {Y}+R^{-1}B\widehat \psi)+C_{1}Z_{1}\right]dt+Z_{1}dW_{1}+Z_{2}dW_{2},\\
Y_T=&\ \zeta.
\end{aligned}
\right.
\end{equation*}
The corresponding Riccati equations are
\begin{equation}\label{Riccati APP2}
\left\{
\begin{aligned}
&\dot\Upsilon-2A\Upsilon-H\Upsilon^2+B^2R^{-1}+C_{1}^2\Upsilon(1+\Upsilon N_{1})^{-1} =0,\\
&\Upsilon_T=0,\\
&\dot \Gamma_{1}+2A\Gamma_{1}-H=0,\\
&\Gamma_{10}=G,\\
&\dot \Gamma_{2}+2A\Gamma_{2}+[B^2R^{-1}+C_{1}^2\Upsilon(1+\Upsilon N_{1})^{-1}]\Gamma_{2}^2-H=0,\\
&\Gamma_{20}=G.
\end{aligned}
\right.
\end{equation}
Equations \eqref{BSDE partial} and \eqref{SDE partial} are reduced to
\begin{equation}\label{BSDE APP2}
\left\{
\begin{aligned}
d\varphi=&\Big[A\varphi+\Upsilon H\widehat \varphi+C_{1}(\eta_{1}-\widehat\eta_{1})+C_{1}(1+\Upsilon N_{1})^{-1}\widehat\eta_{1}\Big]dt+\eta_{1}dW_{1}+\eta_{2}dW_{2},\\
\varphi_T=&\ \zeta,
\end{aligned}
\right.
\end{equation}
and
\begin{equation*}\label{SDE APP2}
\left\{
\begin{aligned}
d\psi=&-\Big\{A\psi+\Gamma_{2}B^2R^{-1}\widehat\psi+\Gamma_{1}C_{1}(\eta_{1}-\widehat\eta_{1})+\Gamma_{2}C_{1}(1+\Upsilon N_{1})^{-1}\left[\widehat \eta_{1}+\Upsilon C_{1}\widehat\psi\right]\Big\}dt
\\&+\Big[(N_{1}-\Gamma_{1})(\eta_{1}-\widehat \eta_{1})+(N_{1}-\Gamma_{2})(I+\Upsilon N_{1})^{-1}\widehat \eta_{1}-C_{1}(1+\Upsilon N_{1})^{-1}(\Gamma_{2}\widehat\varphi+\widehat\psi)\\&-C_{1}\Gamma_{1}(\varphi-\widehat \varphi)-C_{1}(\psi-\widehat \psi)\Big]dW_{1}+(N_{2}-\Gamma_{1})\eta_{2}dW_{2},\\
\psi(0)=&\ 0,
\end{aligned}
\right.
\end{equation*}
respectively.

\begin{figure}[htp]
\begin{center}
\includegraphics[height=6cm,width=10cm]{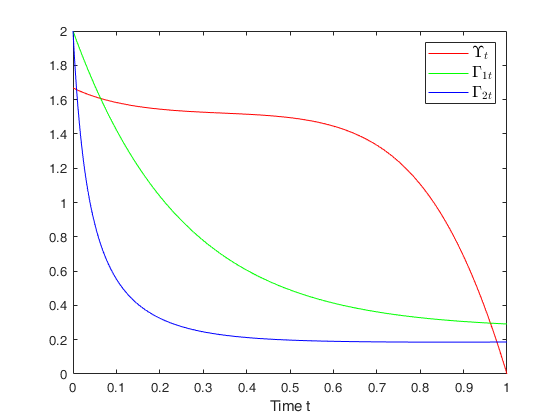}
\caption{The solutions of $\Upsilon, \Gamma_{1}, \Gamma_{2}$}\label{1}
\end{center}\end{figure}

Note that it is hard to obtain  a more explicit expression of $v$ due to the complexity of      \eqref{Riccati APP2} and   \eqref{BSDE APP2}. In the following, we hope to give  numerical solutions for this case with certain particular coefficients.  Let $T=1, A=2, B=3t+2, C_{1}=t-2, G=2, H=e^{-0.05t}, R=2t+1, N_{1}=t(T-t), N_{2}=2$ and $\zeta=T+sin(W_{1T})+cos(2W_{2T})$. Applying Runge-Kutta method, we generate the dynamic simulations
of $\Upsilon, \Gamma_{1}$ and $\Gamma_{2}$, shown in Figure 1.
\begin{figure}[htp]
\begin{center}
\includegraphics[height=6cm,width=10cm]{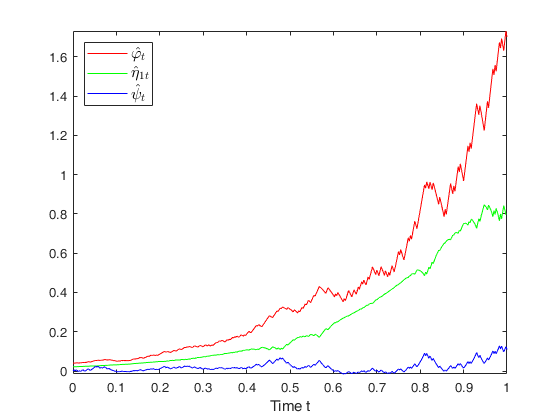}
\caption{Numerical simulations  of $\widehat \varphi$, $\widehat \eta_{1}$ and $\widehat \psi$}\label{2}
\end{center}
\end{figure}

\begin{figure}[htp]
\begin{center}
\includegraphics[height=6cm,width=10cm]{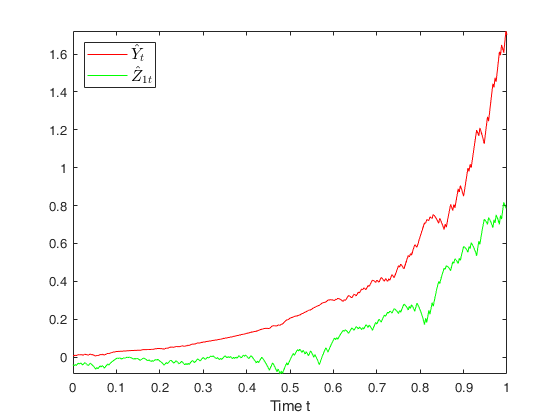}
\caption{Numerical simulations  of $\widehat Y$ and
$\widehat Z_{1}$}\label{3}
\end{center}
\end{figure}

\begin{figure}[htp]
\begin{center}
\includegraphics[height=6cm,width=10cm]{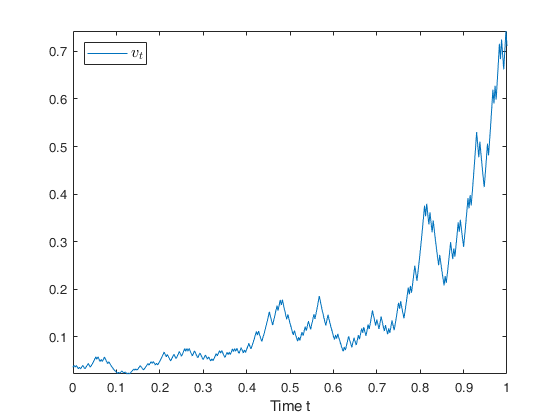}
\caption{Numerical simulation  of $v$}\label{4}
\end{center}
\end{figure}
It seems  that there is no existing literature on   numerical methods of equation \eqref{BSDE  APP2}, which is a BSDE with filtering. Using Theorem 2.1 in Wang et al. \cite{Wangbook} again, we get
\begin{equation*}
\left\{
\begin{aligned}
d\widehat \varphi=&\Big[(A+\Upsilon Q)\widehat \varphi+C_{1}(1+\Upsilon N_{1})^{-1}\widehat\eta_{1}\Big]dt+\widehat\eta_{1}dW_{1},\\
\widehat\varphi_T=&\ \widehat \zeta,\\
d\widehat\psi=&-\Big\{\left[A+\Gamma_{2}B^2R^{-1}+\Gamma_{2}C_{1}^2\Upsilon(1+\Upsilon N_{1})^{-1}\right]\widehat\psi+
\Gamma_{2}C_{1}(1+\Upsilon N_{1})^{-1}\widehat \eta_{1}\Big\}dt
\\&+\Big[(N_{1}-\Gamma_{2})(I+\Upsilon N_{1})^{-1}\widehat \eta_{1}-C_{1}(I+\Upsilon N_{1})^{-1}(\Gamma_{2}\widehat \varphi+\widehat\psi)\Big]dW_{1},\\
\widehat\psi_0=&\ 0.
\end{aligned}
\right.
\end{equation*}
Applying the  numerical method introduced in  Ma et al. \cite{Jin2002Numerical}, we   generate the dynamic simulations of $\widehat \varphi$ and $\widehat \eta_{1}$, shown in Fig. 2. For more information about numerical methods for BSDEs, please refer to Peng and Xu \cite{Peng2011Numerical},
Zhao et al. \cite{Zhao2006A} and the references therein. The simulation of $\widehat \psi$ is also shown in Figure 2.

From Theorem \ref{Relation} and Theorem \ref{main partial}, we have  $\widehat Y=\Upsilon\widehat X+\widehat \varphi$, $\widehat X=-\Gamma_{2}\widehat Y-\widehat\psi$ and $\widehat Z_{1}=(I+\Upsilon N_{1})^{-1}(\widehat\eta_{1}-\Upsilon C_{1}\widehat X)$.  Then the dynamic simulations of $\widehat Y$ and
$\widehat Z_{1}$ are similarly generated, shown in Figure 3.  Further,  from Theorem \ref{main partial}, we   also generate the dynamic simulation of $v$, which is presented in Figure 4.

\section{Conclusion}
 We investigate  an  LQ   control problem of BSDE with partial information, where both the generator of  dynamic system and the cost functional contain diffusion terms $Z_{1}$ and $Z_{2}$. This problem is solved completely and explicitly under some standard conditions. An feedback representation of  optimal control and an  explicit formula of  corresponding optimal cost are given in terms of three Riccati equations, a BSDE with filtering and an  SDE with filtering.  Moreover,  we work out two special scalar-valued control problems to illustrate our theoretical results.

Note that the coefficients in the generator of state equation and the weighting matrices in the cost functional are deterministic. If the coefficients are random, there will be an essential difficulty in solving the case. Since $\mathbb E[A_tY_t|\mathcal F_t^{W_{1}}]=A_t\mathbb E[Y_t|\mathcal F_t^{W_{1}}]$   is no longer true if $A$ is an $\{\mathcal F_t\}_{t\geq 0}$-adapted stochastic process. We will investigate the stochastic case in future.


\begin{thebibliography}{00}

\bibitem{J2000An}
J. M. Bismut, An introductory approach to duality in optimal stochastic control, SIAM Rev.   20 (1978) 62-78.

\bibitem{Pardoux1990Adapted}
E. Pardoux, S.  Peng, Adapted solution of a backward stochastic differential equation, Syst.  Control Lett. 14 (1) (1990) 55-61.

\bibitem{Karoui1997Backward}
N. El. Karoui, S. Peng, M. C. Quenez, Backward stochastic differential equations in finance, Math. Finance 7 (1997)  1-71.

\bibitem{Ma2007Forward}
J. Ma,  J. Yong, Forward-backward stochastic differential equations and their applications, Lecture Notes in Math, Springer-Verlag, New York, 1999.


\bibitem{Kohlmann2000Relationship}
M. Kohlmann, X. Zhou, Relationship between backward stochastic differential equations and stochastic controls: a linear-quadratic approach, SIAM J.   Control   Optim. 38 (5) (2000) 1392-1407.

\bibitem{Peng1993Backward}
S.  Peng, Backward stochastic differential equations and applications to optimal control, Appl. Math. Optim. 27 (2) (1993)  125-144.

\bibitem{Nikolai1999Stochastic}
N. G. Dokuchaev,  X. Zhou, Stochastic controls with terminal contingent conditions, J.  Math. Anal. Appl.  238 (1999)  143-165.

\bibitem{Lim2001Linear}
A. E. B. Lim, X.  Zhou, Linear-quadratic control of backward stochastic differential equations, SIAM J. Control  Optim. 40 (2)  (2001) 450-474.

\bibitem{Li2017Linear}
X. Li, J. Sun,  J. Xiong, Linear quadratic optimal control problems for mean-field backward stochastic differential equations,  Appl. Math. Optim. 80 (2019) 223-250.

\bibitem{Huang2016Backward}
J. Huang, S. Wang, Z. Wu, Backward mean-field Linear-quadratic-gaussian (LQG) games: full and partial information,  IEEE Trans. Automat. Control  60 (12)  (2016) 3784-3796.

\bibitem{Kai2018Linear}
K. Du, J. Huang, Z. Wu, Linear quadratic mean-field-game of backward stochastic differential systems,  Math. Control Relat. Fields  8 (2018)  653-678.

\bibitem{Du2019Linear}
K. Du, Z. Wu,  Linear-quadratic stackelberg game for mean-field backward stochastic differential system and application,  Math. Probl. Eng. 17 (2019)  1-17.



\bibitem{Hu2008Partial}
Y.  Hu, B. \O ksendal, Partial information linear quadratic control for jump diffusions, SIAM J. Control  Optim. 47 (4) (2008) 1744-1761.

\bibitem{Zhen2010A}
Z. Wu, A maximum principle for partially observed optimal control of forward-backward stochastic control systems,  Sci. China Infor. Sci.  53 (11)  (2010) 2205-2214.

\bibitem{huang2009a}
J.  Huang, G.  Wang, J. Xiong, A maximum principle for partial information backward stochastic control problems
with applications,  SIAM J.   Control  Optim. 48 (4)  (2009) 2106-2117.


\bibitem{Wang2015A}
G.  Wang, Z. Wu, J. Xiong, A linear-quadratic optimal control problem of forward-backward stochastic differential equations with partial information, IEEE Trans. Automat. Control 60 (11)  (2015) 2904-2916.

\bibitem{Wang2017An}
G.  Wang, H. Xiao, G.  Xing, An optimal control problem for mean-field forward-backward stochastic differential equation with noisy observation,  Automatica  86 (2017) 104-109.

\bibitem{wang2018a}
G.  Wang, H. Xiao,  J. Xiong, A kind of LQ non-zero sum differential game of backward stochastic differential equation with asymmetric information,  Automatica 97 (2018)  346-352.



\bibitem{Yong1999}
J.  Yong,  X. Zhou, Stochastic Controls: Hamiltonian Systems and HJB Equations, Springer-Verlag, New York, 1999.



\bibitem{Wangbook}
G.  Wang, Z. Wu,  J. Xiong,  An Introduction to Optimal Control of FBSDE with Incomplete Information, Springer-Verlag, New York, 2018.


\bibitem{Xiong2008An}
J. Xiong, An Introduction to Stochastic Filtering Theory,  Oxford University Press, London, 2008.


\bibitem{Liptser}
R. S. Liptser, A. N. Shiryayev, Statistics of Random Processes,
Springer-Verlag,  New York, 1977.





\bibitem{2020A}
P.  Huang, G. Wang,  H.  Zhang,  A partial information linear-quadratic optimal control problem of backward stochastic differential equation with its applications,  Sci. China Infor. Sci.  63 (9) (2020) 1-14.


\bibitem{Jin2002Numerical}
J. Ma, P.  Protter, J. S. Martin,  S.  Torres, Numerical method for backward stochastic differential equations,  Ann. Appl. Probab.  12 (2002)  302-316.





\bibitem{Peng2011Numerical}
S.  Peng, M.  Xu, Numerical algorithms for backward stochastic differential equations with 1-d brownian motion: Convergence and simulations,  ESAIM  Math. Model.  Numer. Anal.  45 (2011)  335-360.


\bibitem{Zhao2006A}
W. Zhao, L.  Chen,  S. Peng, A new kind of accurate numerical method for backward stochastic differential equations,  SIAM J. Sci. Comput. 28 (4)  (2006) 1563-1581.

\end{thebibliography}
\end{document}